# Newton-Mandelbrot set and Murase-Mandelbrot set

Shunji Horiguchi

Abstract. The Wasan(Japanese mathematics in the Edo period(1603-1868: national isolation)) of 1673 leads to three extensions of Mandelbrot's recurrence formula and four extensions of Newton's method. Furthermore, two extended Newton's methods relate to one of the extended Mandelbrot's recurrence formulas(Theorem 3.1). We lead three types of extended Mandelbrot sets: the extended Newton-Mandelbrot set, the Newton-Mandelbrot set, and the Murase-Mandelbrot set. Next, we show that the complex plane can be continuously transformed into a Mandelbrot set or a Multibrot set using the Newton-Mandelbrot set(Theorem 4.4). Next, the Murase-Mandelbrot set becomes a connected, closed set(Proposition 4.12). Finally, examples of Newton-Mandelbrot sets and Murase-Mandelbrot sets are given in §5. These show the originality of Wasan.

## 1. Extended Murase-Mandelbrot recurrence formulas 1,2, and 3.

In 1673, Murase Yoshimasu made a hearth(Fig. 1) with 14 lengths from four rectangular solids. From this hearth volume 192, he created two recurrence formulas(the first and second methods) for the square $x^2$ of the thickness $x$ of the rectangular solid and calculated the square root $x$ by abacus([**1**]).

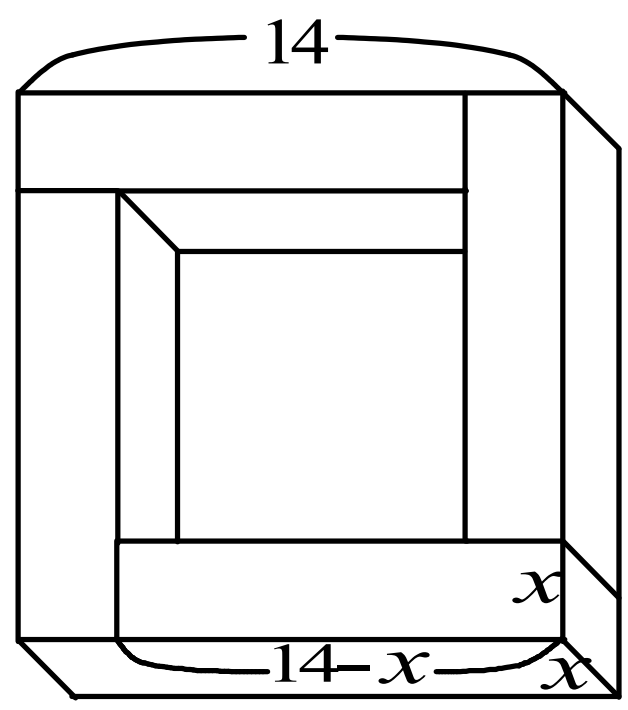

Fig. 1. Hearth

The equation for the volume of the hearth is this.

$$4\cdot x^2(14(1)-x)=4\cdot 48(c) \qquad (1.1)$$

For simplicity, replace 14 by 1 and 48 by $c$ to obtain the following equation.

$$x^3-x^2+c=0 \qquad (1.2)$$

The first method $\qquad x_{k+1}^2 = x_k^3 + c,\ x_0 = 0, k = 0,1,2,\cdots \quad (1.3)$

$$(\Leftrightarrow x_{k+1} = (x_k^3 + c)^{1/2})$$

The second method $\qquad x_{k+1}^2 = \dfrac{c}{1-x_k} \quad (1.4)$

The third method is an equation (1.5), and we obtain the recurrence formula (1.6) from this.

$$c - x^3 = (1 - 2x)x^2 \qquad (1.5)$$

$$x_{k+1}^2 = \frac{c - x_k^3}{1 - 2x_k} \qquad (1.6)$$

Murase changed the third method into the equation (1.5). Formulas (1.3),(1.4), and (1.5) lead to the Horner's method([**4**,**5**]).

**Proposition 1.1.** *Expanded recurrence formula for the volume of the hearth. Formula* (1.7) *is the recurrence formula that the area of a rectangle, the volume of a rectangular solid, and the volume of a* $p$+1*-dimensional rectangle solid, in the order* $p$=1,2, *and* $p$(≥3, *natural number*).

$$x_{k+1}^p = x_k^{p+1} + c \qquad (1.7)$$

**Proposition 1.2.** *Formula* (1.8) *is an extended recurrence formula of formulas* (1.3), (1.4), *and* (1.6) *where m is a real number.*

$$x_{k+1}^{2} = \frac{c + (1-m)x_k^3}{1 - m x_k}, m \in \mathbb{R} \quad (1.8)$$

This means that reducing the length of a side by $x$ reduces its volume by $x^3$.

Replacing the real variable $x$ in (1.8) by the complex variable $z$, we obtain the formula (1.9). This becomes $z_{k+1}^2 = z_k^3 + c$ when $m$=0. This can be regarded as an extension of the Mandelbrot recurrence formula $z_{k+1} = z_k^2 + c$(see def.1.3).

$$z_{k+1}^{2} = \frac{c + (1-m)z_k^3}{1 - m z_k} \quad (1.9)$$

Hence, we go to the definition.

**Definition 1.3.** Let $p$ be a natural number.

$$f(z) = z^{p+1} - z^p + c \quad (1.10)$$

is called the extended Murase-Mandelbrot function.

Let $m$ be a real number.

$$z_{k+1}^{p} = \frac{c + (1-m)z_k^{p+1}}{1 - m z_k} = \frac{(m-1)z_k^{p+1} - c}{m z_k - 1} \quad (1.11)$$

is called the extended Murase-Mandelbrot recurrence formula 1.

$$z_{k+1}^{p} = z_k^{p+1} + c \Leftrightarrow z_{k+1} = (z_k^{p+1} + c)^{1/p} \quad (1.12)$$

is called the extended Murase-Mandelbrot recurrence formula 2. In particular, if $p$=1, then formula (1.12) is called the Mandelbrot recurrence formula.

Let $m$ and $n$ be real numbers that are not 0.

$$z_{k+1} = (z_k^m + c)^n \quad (1.13)$$

is called the extended Murase-Mandelbrot recurrence formula 3.

Assume $x_0$=0 in these (1.11), (1.12), and (1.13).

**Proposition 1.4.** *The domain of* $w=(z^{p+1}+c)^{1/p}(p\geq 2)$ *is a Riemann surface.*

**Proposition 1.5.** *Let* $f_{m,c}(z)=z^m+c$, $z^{1/n}=t$. *Then we obtain*

$$t_{k+1}=f_{m,c}(t_k^{\ n})=(t_k^{\ n})^m+c,\ i.e.,\ z_{k+1}=(z_k^{\ m}+c)^n.$$

## 2. Four extended Newton's methods.

**Definition 2.1.** *Let*

$$f(x)= a_0x^p + a_1x^{p-1} + \cdots + a_{p-1}x + a_p,\ p(\geq 2) \in \mathbb{N} \quad (2.1)$$

be a $p$-th degree polynomial where $a_0,\cdots, a_p$ are real numbers. Let the $j$-th term of $i$-th derived function $f^{(i)}(x)(i=1,\cdots,p)$ of $f(x)$ be $b_j x^{p-j-i}(j=0,1,2,\cdots,p-i)$). Replace the coefficient $b_j$ with a formula including a real variable m or constant and denote such $b_j$ with $b_j(m)$. Next, define such $f^{(i)}(x)$ in $[f^{(i)}(x)]_m$.

For understanding $[f^{(i)}(x)]_m$, go to the next definition and Theorem 2.5.

**Definition 2.2.** We give four extended Newton's methods from Murase's formulas. We choose an initial value $x_o$ near the root of $f(x)$=0.
(1) The following formula to approximate the solution of $f(x)$=0 is called the extended Newton's method 1. Here, $q$ and $\lambda$ are real numbers that are not 0.

$$x_{k+1}^q = x_k^q - \lambda x_k^r \frac{f(x_k)}{f^{(i)}(x_k)},\ q,\lambda(\neq 0) \in \mathbb{R},\ r \in \mathbb{R},\ (k = 0,1,2,\cdots) \quad (2.2)$$

In particular if $q=\lambda=i=1$, $r=0$ then it becomes the Newton(circa 1669.)-Rapson (circa 1690.) method.
(2) Let $f(x)$ be the $p$-th degree polynomial (2.1). The following formula is called the extended Newton's method 2. Here, $q$ and $\lambda$ are real numbers that are not 0.

$$x_{k+1}^{q} = x_k^{q} - \lambda x_k^{r} \frac{f(x_k)}{[f^{(i)}(x_k)]_m}, \; q, \lambda (\neq 0) \in \mathbb{R}, m, r \in \mathbb{R} \quad (2.3)$$

In particular if $q=\lambda=i=1$, $r=0$, $b_j(m)=(p-j)a_j x^{p-(j+1)}(j=0,\cdots,p-1)$ then it becomes the Newton-Rapson method.((2.3) is related to the Horner's method([**4**]).)
(3) The following formula is called the extended Newton's method 3(Tsuchikura-Horiguchi's method).

$$x_{k+1}^{q} = x_k^{q} - q x_k^{q-1} \frac{f(x_k)}{f'(x_k)}, \; q(\neq 0) \in \mathbb{R} \quad (2.4)([\mathbf{3}])$$

(2.2) and (2.3) become Tsuchikura-Horiguchi's method in the special cases.
(4) The following formula is called the extended Newton's method 4(Murase-Newton's method).

$$x_{k+1} = x_k - \lambda \frac{f(x_k)}{f'(x_k)}, \; \lambda(\neq 0) \in \mathbb{R} \quad (2.5)$$

If $\lambda(>1)$ is a natural number, then (2.5) is called Schröder's method([**11**]).

Note that the extended Newton's methods are also valid for complex variables.

**Proposition 2.3.** *If $q(\geq 2)$ is a natural number, then formula* (2.4) *is the first term and the second term of the binomial expansion of the q-square in Newton's method*([**2**]).

**Proposition 2.4.** *From $f(x) = 0$ in formula* (2.1)*, we obtain the following recurrence formula.*

$$x_{k+1}^{p-1} = \frac{(m-a_0)x_k^{p} - a_2 x_k^{p-2} - a_3 x_k^{p-3} - \cdots - a_{p-1} x_k - a_p}{m x_k + a_1}, m \in \mathbb{R} \quad (2.6)$$

**Theorem 2.5.** *The formula* (2.7) *becomes formula* (2.6).

$$x_{k+1}^{p-1} = x_k^{p-1} - (p-1)! \frac{a_0 x_k^p + a_1 x_k^{p-1} + \cdots + a_{p-1} x_k + a_p}{[f^{(p-1)}(x_k)]_m} \quad (2.7)$$

*Here,* $[f^{(p-1)}(x_k)]_m := (p-1)! m x_k + (p-1)! a_1$.

*In particular, formula* (2.7) *becomes Newton's method when* $f(x)=a_0x^2+a_1x+a_2$, $m=2a_0$.

Proof. By replacing $f^{(p-1)}(x)=p!a_0x+(p-1)!a_1$ with $[f^{(p-1)}(x)]_m$, (2.6) is obtained. □

**Theorem 2.6([3]).** *The extended Newton's method* 3 (2.4) *has quadratic convergence when the root of* $f(x)=0$ *is a single root and linear convergence when it is a multiple root.*

The speeds of convergence of Newton's Method and the extended Newton's method 3 can be compared using the convex-concave and curvature of the curve([**3**]).

**Proposition 2.7.** *The following formula is derived from the extended Newton's method* 4 (2.5).

$$\frac{\lambda-1}{\lambda} x_k + \frac{1}{\lambda} x_{k+1} = x_k - \frac{f(x_k)}{f'(x_k)} = {}^N x_{k+1} (\text{denoted as } {}^N x_{k+1}) \quad (2.8)$$

*If* $\lambda$ *is* $1<\lambda$, *then this formula means that Newton's method* ${}^N x_{k+1}$ *is the interior division point that divides the line segment* $x_k x_{k+1}$ *into* $1:\lambda-1$.

*If* $\lambda$ *is* $0<\lambda<1$, *then this formula means that Newton's method* ${}^N x_{k+1}$ *is the exterior division point that divides the line segment* $x_k x_{k+1}$ *into* $1/(1-\lambda):1$.

## 3. Relations between the extended Murase-Mandelbrot recurrence formula 1 and the extended Newton's methods 2 and 3.

**Theorem 3.1.** *Let* $f(z)=z^{p+1}-z^p+c$, $p \in \mathbb{N}$.

(1) *Applying* $f(z)$ *to the extended Newton's method* (2.4), *the obtained formula can be extended to the extended Murase-Mandelbrot recurrence formula* 1 (1.11).

(2) *Applying* $f(z)$ *to the extended Newton's method* (2.7)*, the obtained formula is the extended Murase-Mandelbrot recurrence formula* 1 (1.11).

Proof. (1) Applying $f(z)$ to the extended Newton's method (2.4), we obtain the next formula.

$$z_{k+1}^p = z_k^p - pz_k^{p-1}\frac{z_k^{p+1} - z_k^p + c}{(p+1)z_k^p - pz_k^{p-1}}$$

$$= \frac{z_k^{p+1} - pc}{(p+1)z_k - p} = \frac{(1/p)z_k^{p+1} - c}{(1+1/p)z_k - 1}$$

Here, by replacing $1+1/p$ with any real number $m$ and $1/p$ with $m-1$, we obtain the formula (1.11).

(2) Let the extended Newton's method (2.7) be

$$z_{k+1}^p = z_k^p - p!\frac{z_k^{p+1} - z_k^p + c}{[f^{(p)}(z_k)]_m},$$

here $$[f^{(p)}(z_k)]_m := p!(mz_k - 1).$$

By computing this formula, we obtain formula 1 (1.11). □

## 4. Newton-Mandelbrot set and Murase-Mandelbrot set.

The following are the extended Newton-Mandelbrot set and the Newton-Mandelbrot set.

**Definition 4.1.** Let $\mathbb{C}$ be the complex plane and $n(>1)$ a natural number. We define $M_n$ as follows.

$$M_n := \{c \in \mathbb{C} \mid z_0 = 0,\ z_{k+1} = z_k^n + c,\ \lim_{k\to\infty} |z_k| < \infty\} \quad (4.1)$$

Here, $M_2$ is called the Mandelbrot set(1979), and $M_n(n>2)$ the Multibrot set.

The picture of $M_2$($M_6$ resp.) is Fig. $M_{6,c,1/3}$($M_{3,c,2}$ resp.) in 5.2 of §5.

**Theorem 4.2.** *The Mandelbrot set is a connected, closed set*(*Douady-Hubbard* [**8**])*, and the Multibrot set is also a connected, closed set.*

By Theorem 3.1, the two extended Newton's methods lead to the extended Murase-Mandelbrot recurrence formula 1 (1.11), respectively. Hence, we go to the definition. Here, we give two types of extended Mandelbrot sets.

**Definition 4.3.** Let $p$ be a natural number and $m$ a real number.

$$ENM_{p,m,p+1,c} := \left\{ c \in \mathbb{C} \mid z_0 = 0,\ z_{k+1}^p = \frac{(m-1)z_k^{p+1} - c}{mz_k - 1},\ \lim_{k\to\infty} |z_k| < \infty \right\} \quad (4.2)$$

is called the extended Newton-Mandelbrot set. In particular, $ENM_{1,0,2,c}$ is the Mandelbrot set.

Let $p(>1)$ be a natural number and $m$ a real number.

$$NM_{p,m,c} := \left\{ c \in \mathbb{C} \mid z_0 = 0,\ z_{k+1} = \frac{(m-1)z_k^p - c}{mz_k^{p-1} - 1},\ \lim_{k\to\infty} |z_k| < \infty \right\} \quad (4.3)$$

It is called the Newton-Mandelbrot set($NM$-set). In particular, $NM_{2,0,c}$ is the Mandelbrot set $M_2$ and $NM_{p,0,c}$ $(p>2)$ is the Multibrot set $M_p$.

**Theorem 4.4.** *The complex plane $\mathbb{C}$ can be continuously transformed to Mandelbrot set $M_2$ or Multibrot set $M_p(p>2)$.*

Proof. If we apply $f(z)=z^p-z+c(p>1)$ to Newton's method, then the following is obtained.

$$z_{k+1} = z_k - \frac{z_k^p - z_k + c}{(z_k^p - z_k + c)'} = \frac{(p-1)z_k^p - c}{pz_k^{p-1} - 1} \left( = \frac{(m-1)z_k^p - c}{mz_k^{p-1} - 1} \right)$$

In this case, let $z_0 = 0$ and $z_1 = c$. This is the formula of $NM_{p,p,c}$, and it is the following set.

$$NM_{p,p,c}=\{c(\in \mathbb{C}) \text{ such that } z_k \text{ converges to the root of } f(z)=0 \text{ or } z_k \text{ performs various actions within } \mathbb{C}\}$$

Let $c_{p,k}$ be the result obtained by repeating the initial value $c_p$ $k$ times. Let $c_{p,k}=z_k=(1/p)^{1/(p-1)}$. Then $\mathbb{C}$ is represented as follows.

$$\mathbb{C}=NM_{p,p,c}\cup\{c_{p,k}\}$$

Here, if we continuously move from $p$ to 0, then $c_{0,k}=z_k=(1/p)^{1/(p-1)}=\infty$ ($\Leftrightarrow pz_k^{p-1}-1=0$) and $NM_{p,0,c}=M_p$ is obtained from $\mathbb{C}$. □

**Example 4.5.** If $f(z)=z^3-z+c$, then Newton's method for $f(z)$ is this.

$$z_{k+1}=z_k-\frac{z_k^3-z_k+c}{(z_k^3-z_k+c)'}=\frac{(3-1)z_k^3-c}{3z_k^2-1}\left(=\frac{(m-1)z_k^3-c}{mz_k^2-1}\right)$$

This is the formula of $NM_{3,3,c}$. From $3z_k^2-1=0$, $z_k=\pm(1/3)^{1/2}$ is obtained. Let $c_{3,k}=z_k=\pm(1/3)^{1/2}$. Then $\mathbb{C}$ is represented as follows.

$$\mathbb{C}=NM_{3,3,c}\cup\{c_{3,k}\}$$

Here, if we continuously move from 3 to 0, then $|c_{0,k}|=|z_k|=|\pm(1/3)^{1/2}|=\infty$ and $\mathbb{C}$ transfers to ($NM_{3,0,c}=$)$M_3$ continuously.

We show $NM_{2,m,c}$ approaching $M_2$ in §5.

Next, we will discuss the Fatou set, the Julia set and Newton's map.
Newton's map for $f_{p,c}(z)=z^p-z+c$ is this. Here, add the parameter $m$ to the notation $f_{p,c}$ to write it as $f_{p,m,c}$.

$$Nf_{p,(m=)p,c}(z) = z - \frac{z^p - z + c}{(z^p - z + c)'} = \frac{(p(m)-1)z^p - c}{p(m)z^{p-1} - 1}$$

Let $F(Nf_{p,p,c})$, $J(Nf_{p,p,c})$ be the Fatou set, the Julia set for $Nf_{p,p,c}$ in the Riemann sphere, respectively. Let $c_p=z=(1/p)^{1/(p-1)}$. Then $\mathbb{C}$ is represented as follows.

$$\mathbb{C}=(F(Nf_{p,p,c})\cap \mathbb{C})\cup (J(Nf_{p,p,c})\cap \mathbb{C}),\ c_p\in J(Nf_{p,p,c})\cap \mathbb{C}$$

$$\text{where } \infty\notin (F(Nf_{p,p,c})\cap \mathbb{C})\cup (J(Nf_{p,p,c})\cap \mathbb{C}).$$

Now, we obtain the following.

**Proposition 4.6.** *If* $p\to 0$, *then* $\mathbb{C}=(F(Nf_{p,p,c})\cap \mathbb{C})\cup (J(Nf_{p,p,c})\cap \mathbb{C})$ continuously transformed to

$$\mathbb{C}=(F(Nf_{p,0,c})\cap \mathbb{C})\cup (J(Nf_{p,0,p})\cap \mathbb{C}),$$

$$c_0=\infty\notin \mathbb{C} \text{ where } Nf_{p,0,c}=z^p+c.$$

The following is the Murase-Mandelbrot set.

**Definition 4.7.** Let $m$ and $n$ be real numbers that are not 0.

$$M_{m,c,n} := \{c\in \mathbb{C} \mid z_0 = 0,\ z_{k+1} = (z_k^m + c)^n,\ \lim_{k\to\infty} |z_k| < \infty\} \quad (4.4)$$

is called the Murase-Mandelbrot set. In particular, $M_{2,c,1}$ is the Mandelbrot set. Note that $M_{p,c,1}=NM_{p,0,c}$.

**Definition 4.8.** We give two deformations for the formula $z_{k+1} = (z_k^m + c)^n$.

*Deformation* I $\quad z_{k+1} = (z_k^m + c)^n \Leftrightarrow z_{k+1}^{1/n} = z_k^m + c. \quad (4.5)$

*Deformation* II $\quad z_{k+1} = z_k^{mn} + c \Leftrightarrow z_{k+1}^{1/n} = z_k^m + c = ((z_{k+1}^{1/n})^n)^m + c. \quad (4.6)$

I and II do not move $c$.

**Theorem 4.9.** *In deformation* II*, let* $z_0=0$*. Then the right-hand sides of* $z_{k+1}$ *and* $z_{k+1}^{1/n}$ *are equal.*

Proof. We prove by mathematical induction on $k$. Let $g_{k+1}(c)=z_{k+1}=z_k^{mn}+c$ where $g_0(c)=z_0=0$.

(i) If $l=k=1$ then $g_k(c)=z_k=z_{k-1}^{mn}+c$ and $z_1=c=g_1(c)$, and

$$z_1^{1/n}=z_0^m+c=(z_0^{1/n})^{nm}+c=0^{nm}+c=c=g_1(c)$$

hold.

(ii) If $l=k$ $z_{k+1}=z_k^{mn}+c=g_{k+1}(c)$ and $z_{k+1}^{1/n}=z_k^m+c=g_{k+1}(c)$ hold, then when $l=k+1$, the first recurrence formula is

$$z_{k+2}=z_{k+1}^{mn}+c=(z_k^{mn}+c)^{mn}+c=g_{k+1}(c)^{mn}+c.$$

And the second one is

$$z_{k+2}^{1/n}=z_{k+1}^m+c=((z_{k+1}^{1/n})^n)^m+c=(g_{k+1}(c)^n)^m+c=g_{k+1}(c)^{mn}+c.$$

Hence, we get $z_{k+2}=z_{k+2}^{1/n}$. Furthermore, the convergence of $z_{k+1}=z_k^{mn}+c$ is the same as that of $z_{k+1}^{1/n}=z_k^m+c$. □

We easily derive the following four formulas and their convergences using I, II, and the above Theorem.

**Lemma 4.10.** *If* $z_0=0$*,* $k\longrightarrow\infty$ $z_{k+1}=z_k^{mn}+c=\alpha$ *then the following hold.*

$$z_{k+1}=(z_k^m+c)^n=\alpha^n,\ z_{k+1}=(z_k^n+c)^m=\alpha^m,$$

$$z_{k+1}=(z_k+c)^{mn}=\alpha^{mn} \qquad (4.7)$$

From the above lemma, we obtain the following directly.

**Theorem 4.11.** *The following four Murase-Mandelbrot sets are equal.*

$$M_{mn,c,1}=M_{m,c,n}=M_{n,c,m}=M_{1,c,mn} \qquad (4.8)$$

**Corollary 4.12.** (1) $M_{2,c,1}=M_{m,c,n}=M_{n,c,m}(mn=2)$ (4.9)

(2) $M_{\sqrt{mn},c,\sqrt{mn}} = M_{mn,c,1} = M_{1,c,mn}\ (m>0, n>0)$ (4.10)

By Theorem 4.11, we obtain the next proposition.

**Proposition 4.13.** *Let m,n be the natural numbers such that* $mn \geq 3$. *Then, the Multibrot set* $M_{mn,c,1}$ *is a connected, closed set. And* $M_{m,c,n}=M_{n,c,m}=M_{1,c,mn}(=M_{mn,c,1})$ *is also a connected, closed set.*

The following proposition gives the calculation of the recurrence formula to obtain the picture of $M_{1,c,mn}$.

**Proposition 4.14.** *Calculation of* $z_{k+1}=(z_k+c)^s$. *Let* $s(\neq 0)$ *be a real number, and* $z=x+iy, c=a+ib, X=x+a, Y=y+b, \alpha=X+iY$. *Then the following holds.*

$$z_{k+1}=(z_k+c)^s=\alpha_k^s=|\alpha_k|^s e^{is(\mathrm{Arg}\alpha_k+2n\pi)},\ n=0,\pm1,\pm2,\cdots \qquad (4.11)$$

*In particular, if* $s(\geq 2)$ *is a natural number, then the following holds.*

$$z_{k+1}=(z_k+c)^s=\alpha_k^s=|\alpha_k|^s e^{i(s\mathrm{Arg}\alpha_k)}(=(|\alpha_k|^2 e^{i(2\mathrm{Arg}\alpha_k)})^{s/2}=((z_k+c)^2)^{s/2}) \qquad (4.12)$$

## 5. Examples of Newton-Mandelbrot sets $NM_{2,m,c}$ and Murase-Mandelbrot sets $M_{m,c,n}$.

**5.1. Examples of $NM_{2,m,c}$.** The recurrence formula of $NM_{2,m,c}$ is this.

$$z_{k+1}=\frac{(m-1)z_k^2-c}{mz_k-1}$$

If $m$ continuously approaches 0, then this approaches $z_{k+1}=z_k^2+c$.

***NM*$_{2,-1,c}$** The formula of $NM_{2,-1,c}$ is this.

$$z_{k+1}=\frac{2z_k^2+c}{z_k+1}$$

This picture is Fig. $NM_{2,-1,c}$.

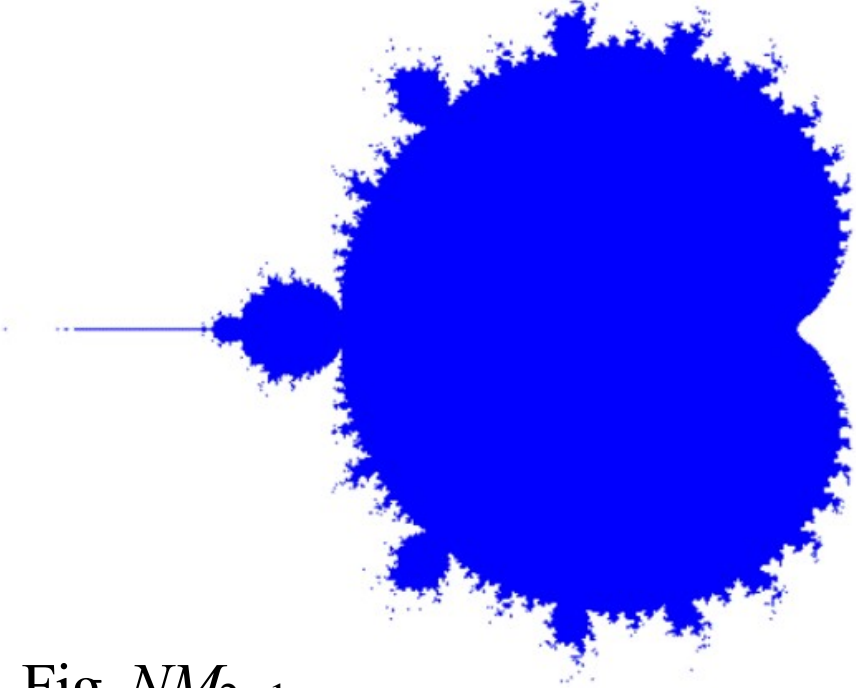

Fig. $NM_{2,-1,c}$.

***NM*$_{2,0.1,c}$** The formula of $NM_{2,0.1,c}$ is this.

$$z_{k+1}=\frac{-0.9z_k^2-c}{0.1z_k-1}$$

This picture is close to the shape of the Mandelbrot set $NM_{2,0,c}(=M_2)$. See Fig. $M_{6,c,1/3}$.

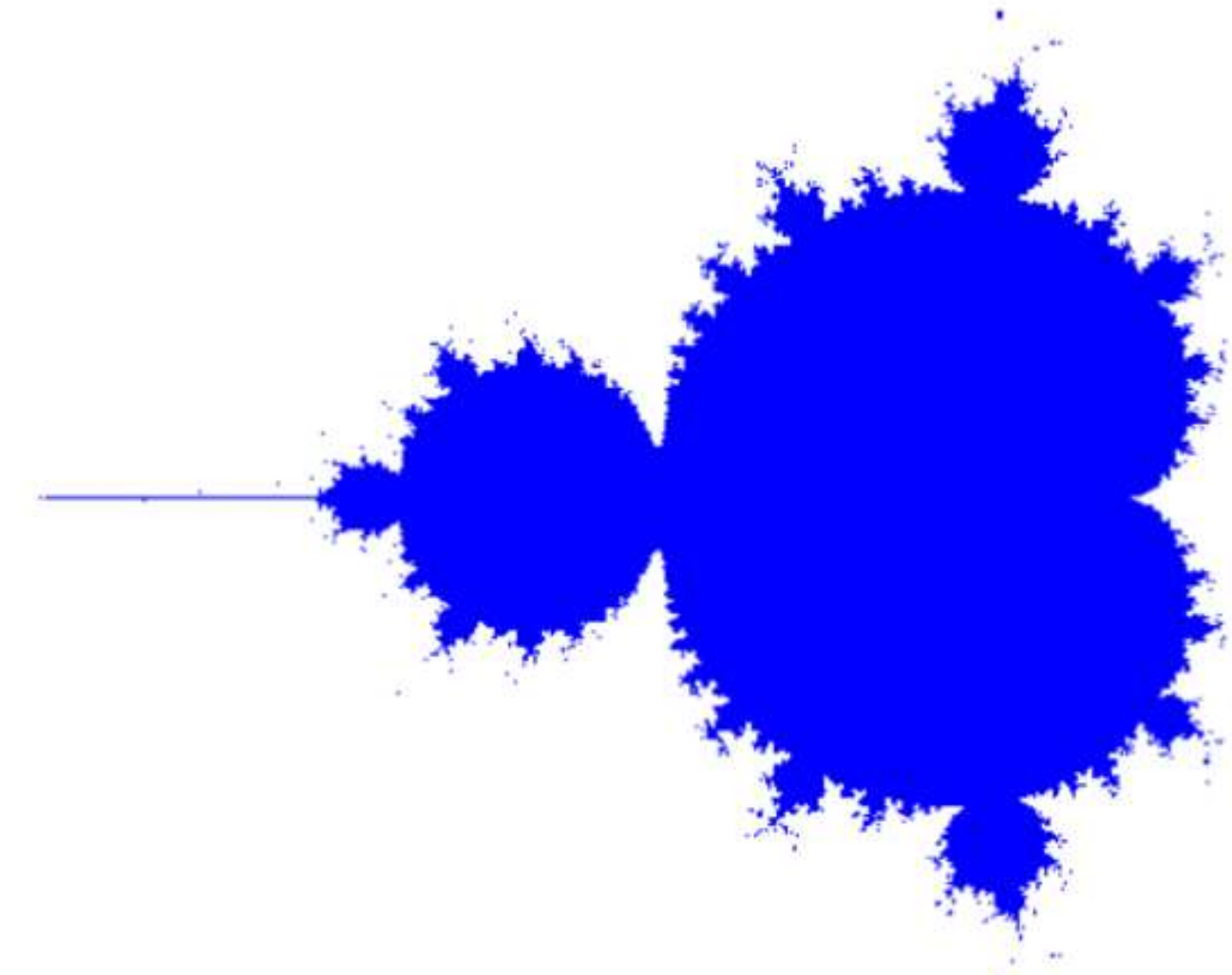

Fig. $NM_{2,0.1,c}$.

***$NM_{2,0.5,c}$*** The formula of $NM_{2,0.5,c}$ is this.

$$z_{k+1} = \frac{-0.5z_k^2 - c}{0.5z_k - 1} = \frac{z_k^2 + 2c}{-z_k + 2}$$

The picture of $NM_{2,0.5,c}$ is this.

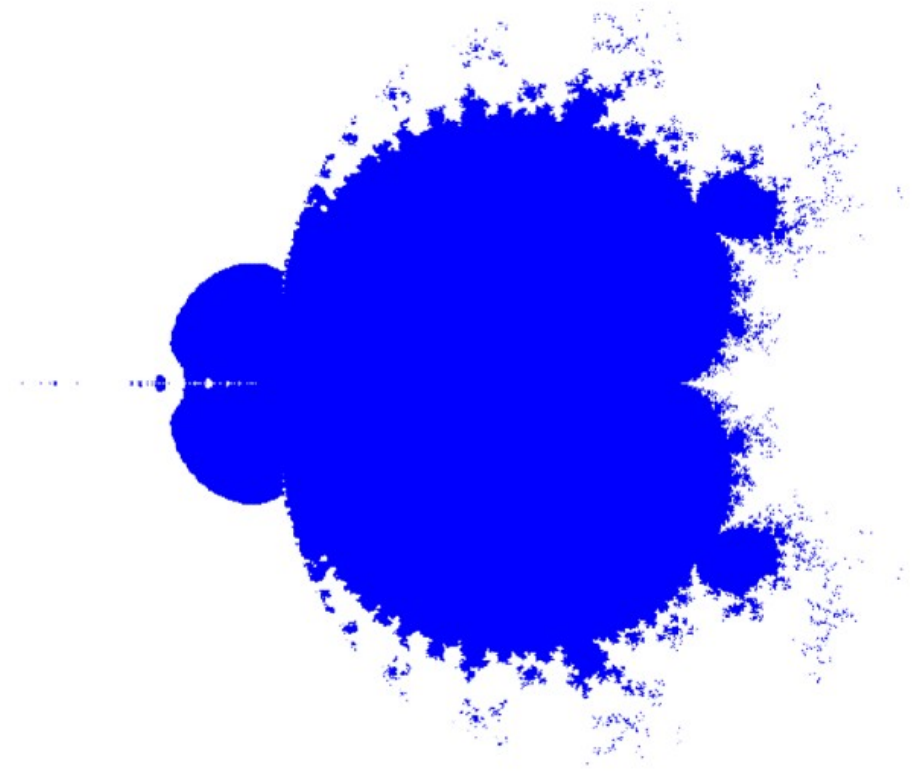

Fig. $NM_{2,0.5,c}$.

***$NM_{2,1,c}$*** The formula of $NM_{2,1,c}$ is this.

$$z_{k+1} = \frac{-c}{z_k - 1}$$

$NM_{2,1,c}$ has a wide convergence range.

**5.2. Examples of $M_{m,c,n}$.** The recurrence formula is $z_{k+1} = (z_k^m + c)^n$.

***$M_{3,c,2}$*** Formula $z_{k+1} = (z_k^3 + c)^2$ deforms into

$$z_{k+1} = z_k^6 + c,\ z_{k+1} = (z_k^2 + c)^3, \text{and}\ \ z_{k+1} = (z_k + c)^6.$$

These four extended Mandelbrot sets are the same. We used the formula

$z_{k+1}=(z_k+c)^6$ to draw this picture. See Fig. $M_{1,c,6}$. The white area is the convergence range. $M_{3,c,2}$ is a connected and closed set.

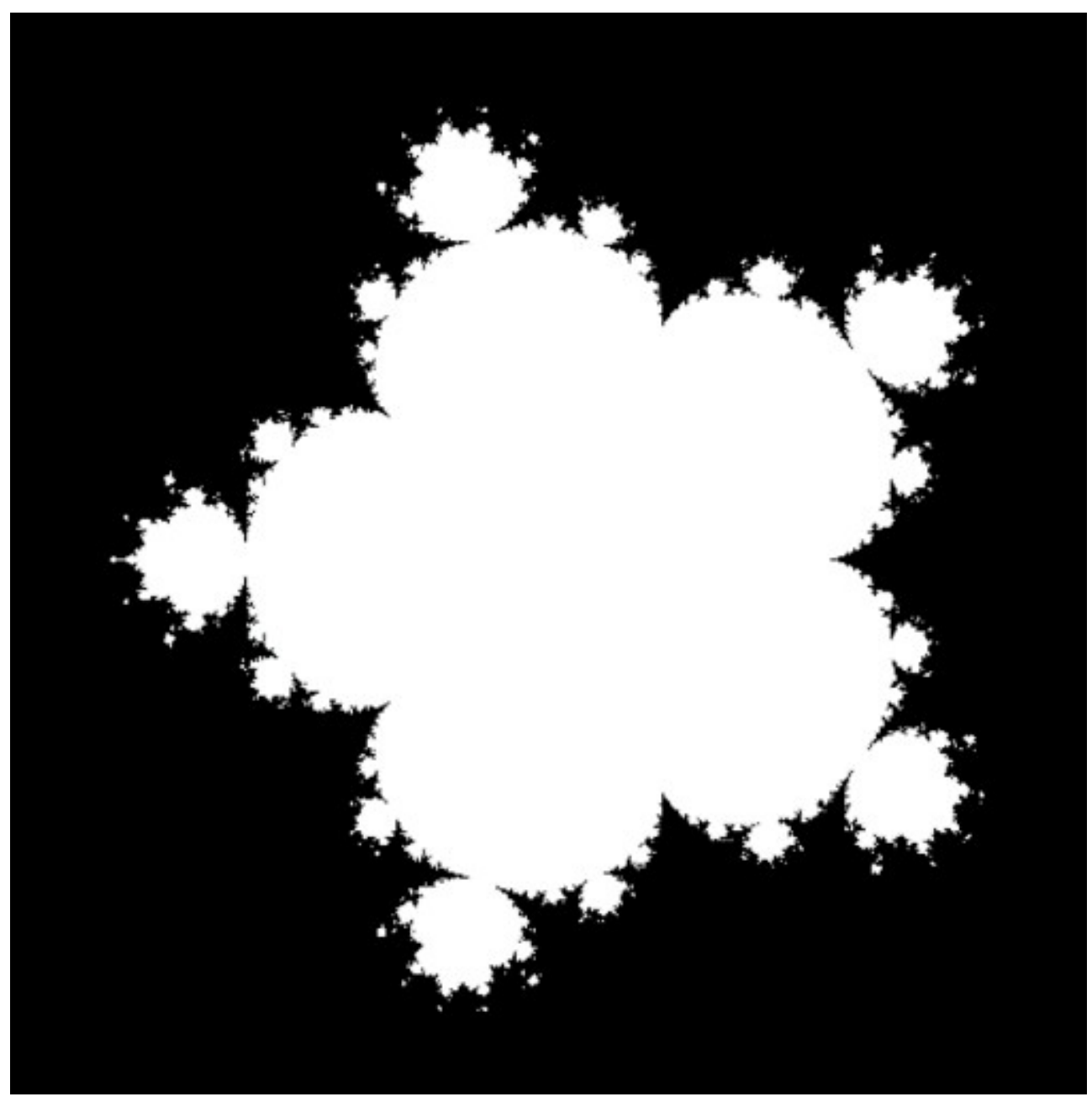

Fig. $M_{1,c,6}(=M_{6,c,1}=M_{2,c,3}=M_{3,c,2})$.

***$M_{6,c,1/3}$*** This is an interesting example. The following deformations are held by deformations Ⅰ and Ⅱ.

$$z_{k+1}=(z_k^6+c)^{1/3} \quad (1) \;\Leftrightarrow\; z_{k+1}=(z_k^{1/3}+c)^6 \quad (2)$$

$$\Leftrightarrow\; z_{k+1}=z_k^2+c \qquad (3) \;\Leftrightarrow\; z_{k+1}=(z_k+c)^2 \quad (4)$$

These formulas give the Mandelbrot set. This picture is Fig. $M_{6,c,1/3}$. This is derived from formula (1). (3) is the formula of the Mandelbrot set. Formula (1) is Murase-Mandelbrot set. The computer calculation in (1) applies formula (4.11)($n$=0,1,2) to $z_k^6+c$. However, the above deformations assert that each picture of these three($n$=0,1,2) is the Mandelbrot set. In fact, when the three

formulas($n$=0,1,2) of the Murase-Mandelbrot sets are computed and drawn by a computer, they become the picture of the Mandelbrot set. If $m$ and $n$ are real numbers such that $mn$=2, then $M_{m,c,n}$ will be $M_{2,c,1}$.

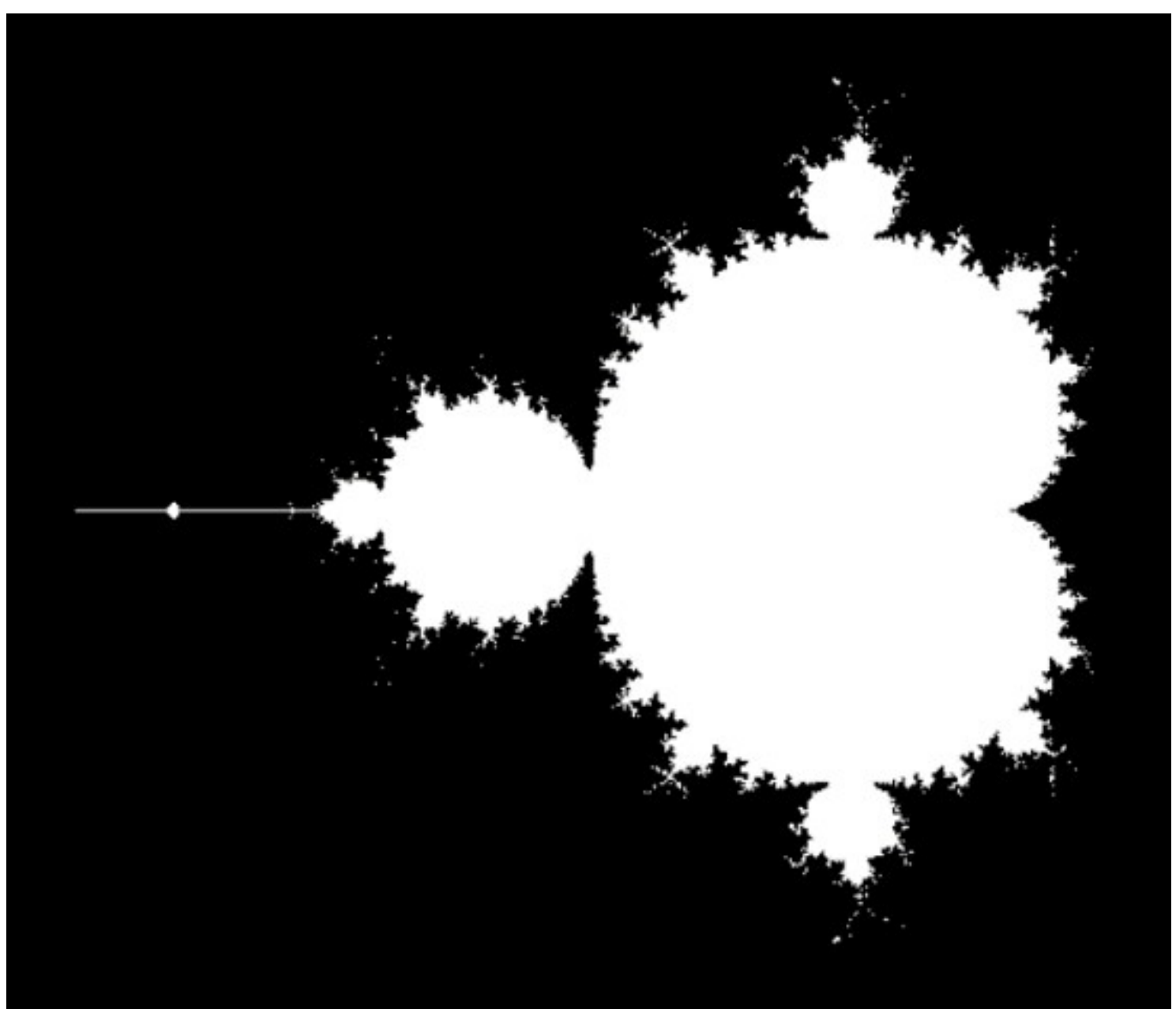

Fig. $M_{6,c,1/3}$. Mandelbrot set $M_{2,c,1}$.

**$M_{3,c,1/2}$** Following deformations hold.

$$z_{k+1} = (z_k^3 + c)^{1/2} \Leftrightarrow z_{k+1} = (z_k^{1/2} + c)^3 \Leftrightarrow z_{k+1} = z_k^{3/2} + c$$

$$\Leftrightarrow z_{k+1} = (z_k + c)^{3/2}$$

Formula (1) gives two pictures of $M_{3,c,1/2}$. The computer calculation is formula (4.11). We obtain two pictures Fig. $M_{3,c,1/2}$,$n$=0,2,4,･･･and Fig. $M_{3,c,1/2}$, $n$=1,3,5,･･･.

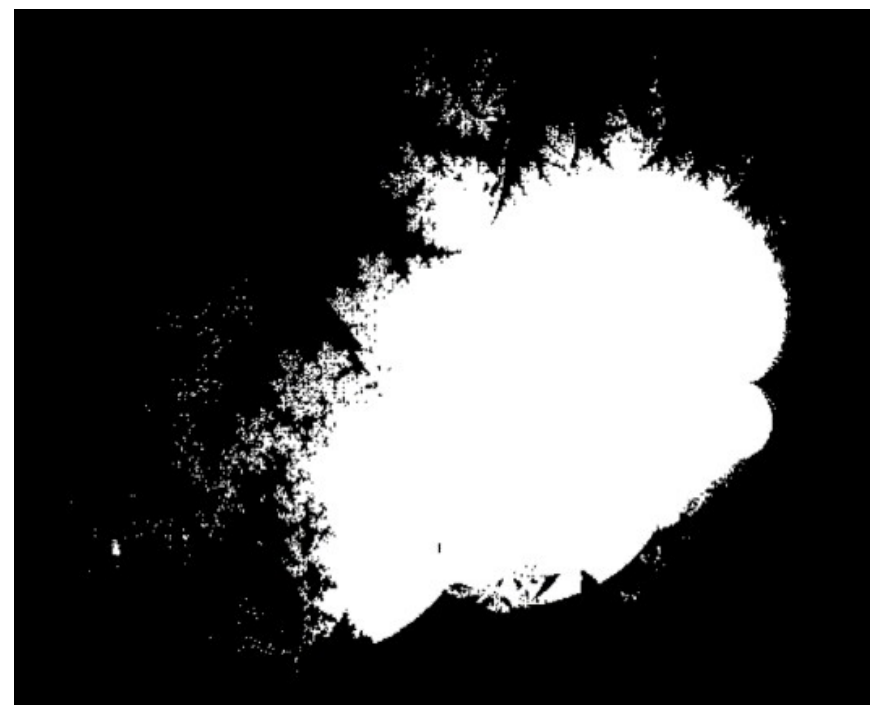

Fig. $M_{3,c,1/2}$, $n$=0,2,4,・・・.

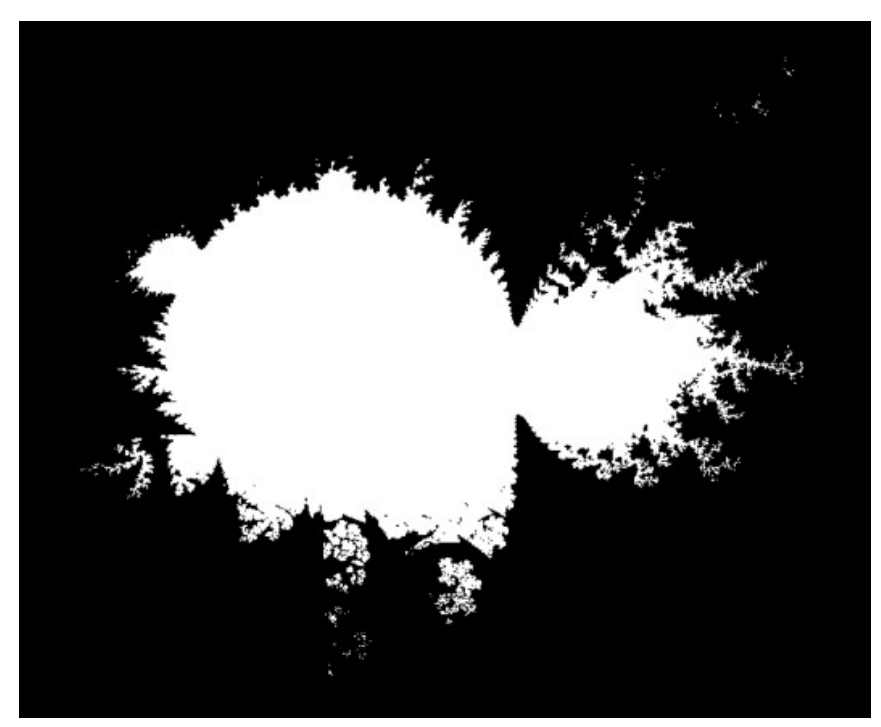

Fig. $M_{3,c,1/2}$, $n$=1,3,5,・・・.

$\boldsymbol{M_{1,c,\sqrt{2}}}$ The following deformation holds.

$$z_{k+1} = (z_k + c)^{\sqrt{2}} \;\; (1) \Leftrightarrow z_{k+1} = z_k^{\sqrt{2}} + c \;\; (2)$$

Formula (1) gives an infinite number of pictures of $M_{1,c,\sqrt{2}}$. The computer calculation of (1) is

$$|\alpha_k|^{\sqrt{2}} e^{i(\sqrt{2}\mathrm{Arg}\,\alpha_k + 2\sqrt{2}n\pi)} \; (n=0,\pm1,\pm2,\cdots).$$

We give four pictures of $M_{1,c,\sqrt{2}}$ when $n$=0,1,2, and 3.

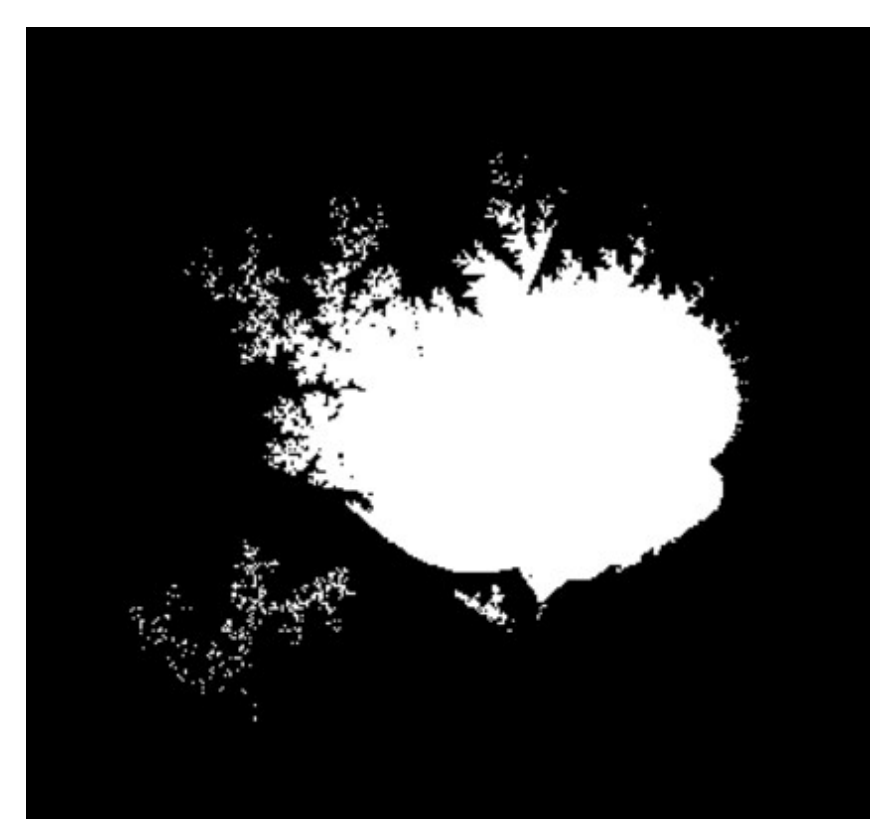

Fig. $M_{1,c,\sqrt{2}}$, $n$=0.

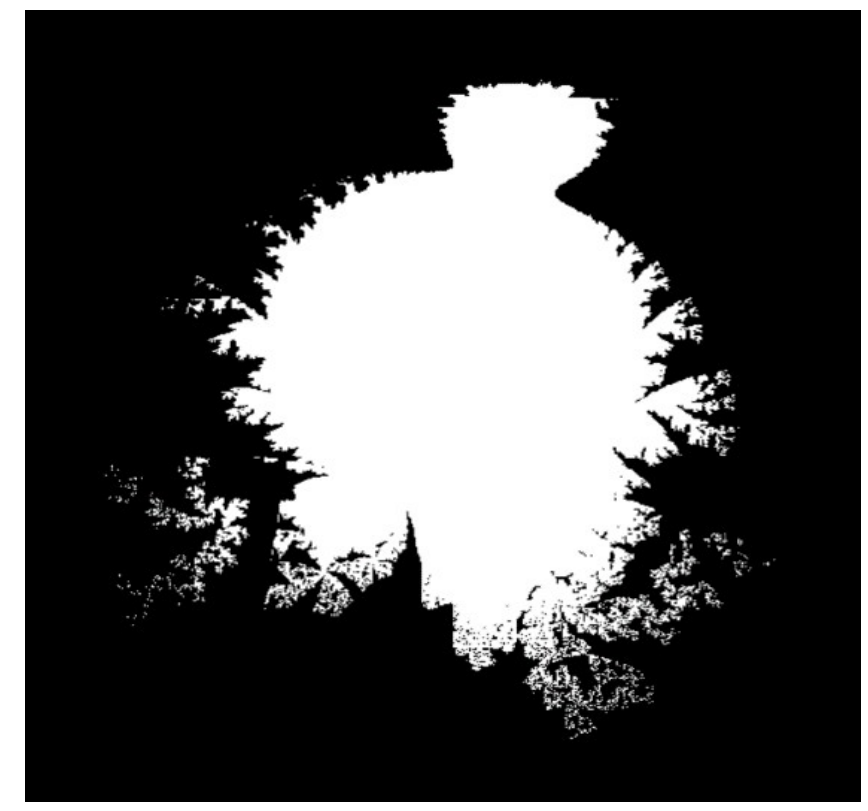

Fig. $M_{1,c,\sqrt{2}}$, $n$=1.

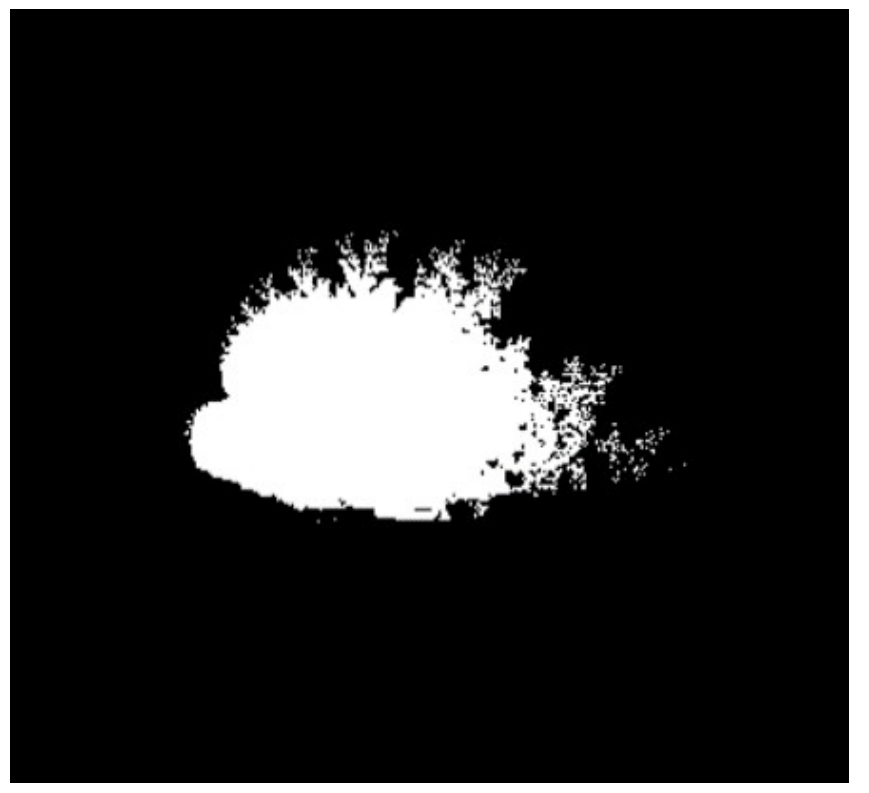

Fig. $M_{1,c,\sqrt{2}}$, $n$=2.

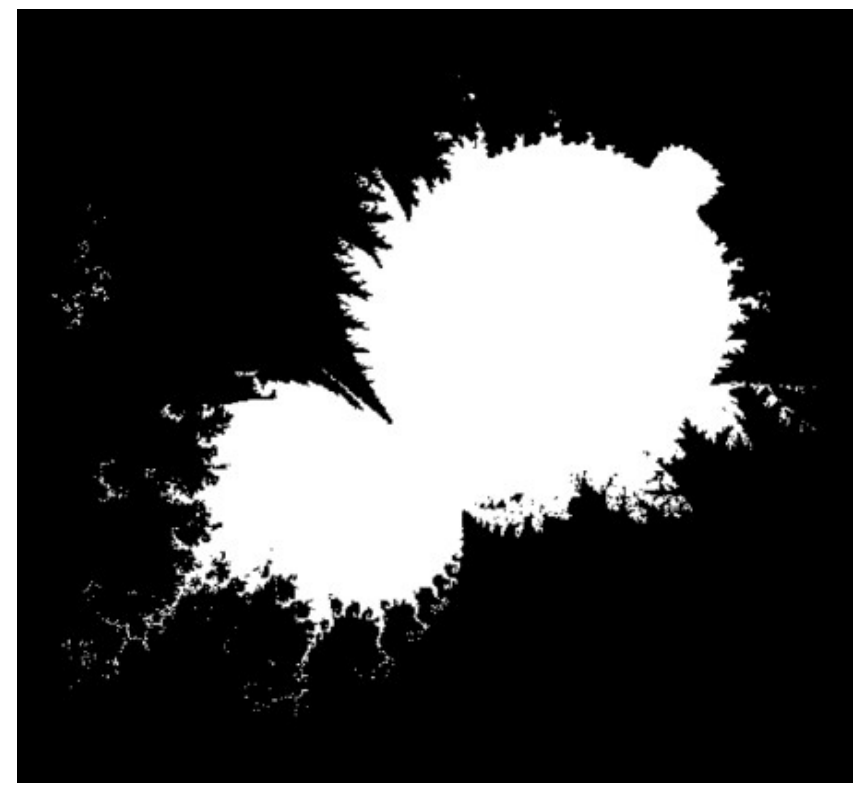

Fig. $M_{1,c,\sqrt{2}}$, $n$=3.

Acknowledgment. I received advice from Professor Mitsuo Morimoto and Professor Tomoki Kawahira. I would like to express my gratitude here.

Shunji HORIGUCHI
JAPAN
shunhori@seagreen.ocn.ne.jp